\documentclass[12pt]{article}

\usepackage{amssymb}

\newtheorem {theorem} {Theorem}

\newtheorem {lemma}  [theorem]{Lemma}

\newcommand{\bbox}{\ \hfill\rule[-1mm]{2mm}{3.2mm}}

\def\div{{\rm div}}

\title{\sc On the stability of limit cycles for planar
differential systems.\thanks{The second author is partially
supported by a MCYT grant number BFM 2002-04236-C02-01 and by a
FPU grant with reference AP2000-3585.}}

\author{{\sc H. Giacomini$^{\ (1)}$ and M.
Grau$^{\ (2)}$}}
\date{}

\begin{document}
\maketitle \vspace{-1.0cm}
\begin{center}
$^{\ (1)}$ Laboratoire de Math\'ematiques et Physique Th\'eorique.
\\ C.N.R.S. UMR 6083. \\ Facult\'e des Sciences et Techniques.
Universit\'e de Tours. \\ Parc de Grandmont, 37200 Tours, FRANCE.
\\ E-mail: {\tt Hector.Giacomini@phys.univ-tours.fr}
\\ $ $ \\
$^{\ (2)}$ Departament de Matem\`atica. Universitat de Lleida. \\
Avda. Jaume II, 69. 25001 Lleida, SPAIN. \\ {\rm E--mail:} {\tt
mtgrau@matematica.udl.es}
\end{center}
\begin{abstract}
We consider a planar differential system $\dot{x}= P(x,y)$,
$\dot{y} = Q(x,y)$, where $P$ and $Q$ are $\mathcal{C}^1$
functions in some open set $\mathcal{U} \subseteq \mathbb{R}^2$,
and $\dot{}=\frac{d}{dt}$. Let $\gamma$ be a periodic orbit of the
system in $\mathcal{U}$. Let $f(x,y): \mathcal{U} \subseteq
\mathbb{R}^2 \to \mathbb{R}$ be a $\mathcal{C}^1$ function such
that
\[ P(x,y) \, \frac{\partial f}{\partial x}(x,y) \, + \, Q(x,y) \,
\frac{\partial f}{\partial y} (x,y) \, = \, k(x,y) \, f(x,y), \]
where $k(x,y)$ is a $\mathcal{C}^1$ function in $\mathcal{U}$ and
$\gamma \subseteq \{ (x,y) \, | \, f(x,y) = 0\}$. We assume that
if $p \in \mathcal{U}$ is such that $f(p)=0$ and $\nabla f(p)=0$,
then $p$ is a singular point.
\par
We prove that $\int_{0}^{T} \left( \frac{\partial P}{\partial x} +
\frac{\partial Q}{\partial y}\right)(\gamma(t)) \, dt= \int_0^{T}
k(\gamma(t)) \, dt$, where $T>0$ is the period of $\gamma$. As an
application, we take profit from this equality to show the
hyperbolicity of the known algebraic limit cycles of quadratic
systems.
\end{abstract}

\section{Introduction \label{sect1}}

We consider a {\em planar differential system}
\begin{equation}
\dot{x} = P(x,y), \quad \dot{y} = Q(x,y), \label{eq1}
\end{equation}
where $P$ and $Q$ are $\mathcal{C}^1$ functions in some open set
$\mathcal{U} \subseteq \mathbb{R}^2$, and $\dot{}=\frac{d}{dt}$. A
singular point of system (\ref{eq1}) is a point $p \in
\mathcal{U}$ such that $P(p)=Q(p)=0$. We assume that all the
singular points of (\ref{eq1}) are isolated.
\par
Given a system (\ref{eq1}), we can always consider its vector
field representation $\mathbf{F}(x,y)=(P(x,y),Q(x,y))$.
\par
We will denote by $\div(x,y)$ the {\em divergence} of system
(\ref{eq1}), that is, $\div = \partial P / \partial x \, + \,
\partial Q / \partial y$.
\newline

We also need to consider the {\em flow} of system (\ref{eq1}),
which we denote by $\Phi_{t}(p)$ and which represents the unique
solution of system (\ref{eq1}) passing through the point $p \in
\mathcal{U} \subseteq \mathbb{R}^2$. We notice that for each $p
\in \mathcal{U}$ there exists an $\epsilon_p >0$ (which may be
$\epsilon_p = + \infty$) such that $t \in (- \epsilon_p,
\epsilon_p)$ is the maximal symmetric interval of existence of the
solution of (\ref{eq1}) passing through $p$. We have that $\frac{d
\Phi_t}{dt} (p) = (P(\Phi_t(p)),Q(\Phi_t(p)))$, for all $p \in
\mathcal{U}$ and $t \in (-\epsilon_p, \epsilon_p)$, and
$\Phi_0(p)=p$. Given $p \in \mathcal{U}$, the function $\Phi(
\cdot, p): (-\epsilon_p, \epsilon_p) \to \mathbb{R}^2$, where
$\Phi(t,p):=\Phi_t(p)$, defines a {\em solution curve} or {\em
orbit} of (\ref{eq1}) through the point $p$.
\par
A {\em limit cycle} of system (\ref{eq1}) is an isolated periodic
orbit. Let $\gamma$ be a limit cycle for system (\ref{eq1}). We
say that $\gamma$ is {\em stable} if there exists a neighborhood
$\mathcal{U}_\gamma \subseteq \mathcal{U}$ of $\gamma$ such that
for all $p \in \mathcal{U}_\gamma$, we have $\lim_{t \to + \infty}
d(\Phi_t(p), \gamma) =0$. As usual, the previous application $d$
is the distance between sets in the Hausdorff sense. Analogously,
we say that $\gamma$ is {\em unstable} if there exists a
neighborhood $\mathcal{U}_\gamma \subseteq \mathcal{U}$ of
$\gamma$ such that for all $p \in \mathcal{U}_\gamma$, we have
$\lim_{t \to - \infty} d(\Phi_t(p), \gamma) =0$.
\par
There might be limit cycles which are neither stable nor unstable.
Using the Jordan curve theorem, which states that any simple
closed curve, as the limit cycle, $\gamma$ separates any
neighborhood $\mathcal{U}_\gamma$ of $\gamma$ into two disjoint
sets having $\gamma$ as a boundary, we can consider
$\mathcal{U}_\gamma$ as the disjoint union of $\mathcal{U}_i \cup
\gamma \cup \mathcal{U}_e$, where $\mathcal{U}_i$ and
$\mathcal{U}_e$ are open sets situated, respectively, in the
interior and exterior of $\gamma$. When for any $p \in
\mathcal{U}_i$ we have $\lim_{t \to + \infty} d(\Phi_t(p), \gamma)
=0$ whereas for any $q \in \mathcal{U}_e$ we have $\lim_{t \to -
\infty} d(\Phi_t(q), \gamma) =0$ (or, the other way round, for any
$p \in \mathcal{U}_i$ we have $\lim_{t \to - \infty} d(\Phi_t(p),
\gamma) =0$ whereas for any $q \in \mathcal{U}_e$ we have $\lim_{t
\to + \infty} d(\Phi_t(q), \gamma) =0$), we say that $\gamma$ is
{\em semi-stable}.
\par
Any limit cycle $\gamma$ of a system (\ref{eq1}) is either stable,
unstable or semi-stable as it is stated in \cite{Perko}. For a
detailed description of the classical known results on limit
cycles see also \cite{Perko}.
\par
The following result, which is stated as a corollary in page 214
of \cite{Perko}, gives a formula to distinguish the stability of a
limit cycle.
\begin{theorem}
Let $\gamma(t)$ be a periodic orbit of system {\rm (\ref{eq1})} of
period $T$. Then, $\gamma$ is a stable limit cycle if \[
\int_{0}^T \div(\gamma(t)) \, dt <0, \] and it is an unstable
limit cycle if \[ \int_{0}^T \div(\gamma(t)) \, dt >0. \] It may
be stable, unstable or semi-stable limit cycle or it may belong to
a continuous band of cycles if this quantity is zero.
\label{thdiv}
\end{theorem}
A sketch of the proof of this theorem is given after the
forthcoming Theorem \ref{thDP}.
\par
When the quantity $\int_{0}^T \div(\gamma(t)) \, dt$ is different
from zero, we say that the limit cycle $\gamma$ is {\em
hyperbolic}.
\par
Since we are considering differential systems (\ref{eq1}) in the
class of functions $\mathcal{C}^1$, we may have a limit cycle
$\gamma$ belonging to a sequence of periodic orbits $\{ \gamma_n
\, , \, n \in \mathbb{N} \}$ with $\gamma_{n+1}$ in the interior
of $\gamma_n$, such that the sequence accumulates to a singular
point, a periodic orbit or a graphic and such that every
trajectory between $\gamma_n$ and $\gamma_{n+1}$ spirals towards
$\gamma_n$ or $\gamma_{n+1}$ as $t \to \pm \infty$. This kind of
phenomena does not exist for analytic systems.
\par
In this work, we give another quantity which equals to $\int_{0}^T
\div(\gamma(t)) \, dt$ for a periodic orbit $\gamma$ defined in an
implicit way, as explained below. This is the main result of the
article and it is stated in Theorem \ref{threl} in the following
section. We can, therefore, distinguish the hyperbolicity of a
limit cycle using two different quantities.
\newline

Given a planar system (\ref{eq1}) (or equivalently its vector
field representation $\mathbf{F}(x,y)=(P(x,y),Q(x,y))$), we define
an {\em invariant curve} as a curve given by $f(x,y)=0$, where $f:
\mathcal{U} \subseteq \mathbb{R}^2 \to \mathbb{R}$ is a
$\mathcal{C}^1$ function in the open set $\mathcal{U}$, non
locally null and such that there exists a $\mathcal{C}^1$ function
in $\mathcal{U}$, denoted by $k(x,y)$, for which
\begin{equation}
P(x,y) \, \frac{\partial f}{\partial x} (x,y) \, + \, Q(x,y) \,
\frac{\partial f}{\partial y}(x,y) \, = \, k(x,y) \, f(x,y),
\label{ci}
\end{equation}
for all $(x,y) \in \mathcal{U}$. The identity (\ref{ci}) can be
rewritten by $\nabla f \cdot \mathbf{F} = k f$. As usual, $\nabla
f$ denotes the gradient vector related to $f(x,y)$, that is,
$\nabla f(x,y)=( \frac{\partial f}{\partial x} (x,y) ,
\frac{\partial f}{\partial y} (x,y) )$, $\mathbf{F}(x,y)$ is the
previously defined vector $(P(x,y),Q(x,y))$, and $\cdot$ denotes
the scalar product. We will denote by $\dot{f}$ or by $\frac{d
f}{d t}$ the function $\nabla f \cdot \mathbf{F}$ once evaluated
on a solution of system (\ref{eq1}).
\par
We will always assume that if $p \in \mathcal{U}$ is such that
$f(p)=0$ and $\nabla f(p)=0$, then $p$ is a singular point of
system (\ref{eq1}). This is a technical hypothesis which
generalizes the notion of not having multiple factors for
algebraic curves. For instance, if we had written that the
periodic orbit $\gamma$ was contained in $f^2(x,y)=0$, then we
would have that $\nabla (f^2) (p)=0$ for all $p \in \gamma$, in
contradiction with the hypothesis.
\par
We notice that, as a particular case, we may have a function
$f(x,y)$ given by a polynomial in $\mathbb{R}[x,y]$. In such a
case, $f(x,y)=0$ is called an {\em invariant algebraic curve}.
When, in addition, the system is polynomial, that is, $P,Q \in
\mathbb{R}[x,y]$, then the function $k(x,y)$ is a real polynomial
called {\em cofactor}. When we consider an algebraic curve, we can
always assume that it is defined by a polynomial $f(x,y)=0$ such
that the decomposition of $f(x,y)$ has no multiple factors. The
same assumption must be done for curves defined by $\mathcal{C}^1$
functions and it is equivalent to the assumption that if $p \in
\mathcal{U}$ is such that $f(p)=0$ and $\nabla f(p)=0$, then $p$
is a singular point of system (\ref{eq1}). More explicitly, assume
that the system (\ref{eq1}) has an invariant algebraic curve given
by $f(x,y)=0$ and assume that the decomposition of $f(x,y)$ in the
ring $\mathbb{R}[x,y]$ has no multiple factors, that is, $f(x,y) =
b_1(x,y) b_2(x,y) \ldots b_k(x,y)$ where $b_j(x,y)$ is an
irreducible polynomial in $\mathbb{R}[x,y]$ and $b_i(x,y) \neq c
b_j(x,y)$ for any $c \in \mathbb{R} - \{0\}$ if $i \neq j$. Let $p
\in \mathcal{U}$ be such that $f(p) = 0$ and $\nabla f(p)=0$.
Since the decomposition of $f(x,y)$ has no multiple factors, we
deduce that $p$ is a singular point of the curve $f(x,y)=0$ and,
hence, it is a singular point of the system (\ref{eq1}).
\par
Our main result, Theorem \ref{threl}, can only be applied when the
periodic orbit $\gamma$ is given in an implicit way, that is, when
there exists an invariant curve $f(x,y)=0$ such that $\gamma
\subseteq \{ (x,y) \, | \, f(x,y) = 0\}$. For instance, let us
consider the following $\mathcal{C}^1$ system defined in all
$\mathbb{R}^2$:
\begin{equation}
\dot{x} =  (x+y) \cos(x) - y (x^2+xy+2y^2), \  \dot{y}  =
 (y-x) (\cos(x) - y^2) + \frac{x^2+y^2}{2} \sin(x),
\label{nalc}
\end{equation}
which has $y^2- \cos(x) = 0$ as invariant curve. We define
$f(x,y):=y^2-\cos(x)$ and we have that $f \in
\mathcal{C}^1(\mathbb{R}^2)$ and that $\nabla
f(x,y)=(\sin(x),2y)$. Therefore, there is no $p \in \mathbb{R}^2$
such that both $f(p)=0$ and $\nabla f(p)=0$. Moreover, $f(x,y)=0$
satisfies equation (\ref{ci}) with $k(x,y)=2y(x-y)-(x+y) \sin(x)$.
The divergence of this system is $\div(x,y)= -4y^2+2\cos(x)-x
\sin(x)$ and $V(x,y)=(x^2+y^2) f(x,y)$ is an inverse integrating
factor. We denote by $\gamma_n$, $n \in \mathbb{Z}$, the oval of
$f(x,y)=0$ belonging to the strip $-\pi/2+ 2 \pi n \leq x \leq
\pi/2 + 2 \pi n$. Each oval $\gamma_n$, with $n \in \mathbb{Z}$,
gives a periodic orbit of (\ref{nalc}) with minimal period
$T_n>0$. The oval $\gamma_0$ is a hyperbolic stable limit cycle
for system (\ref{nalc}), which can be shown just applying Theorem
\ref{thdiv}. We have, after some easy computations, that
\[ \int_{0}^{T_n} \div(\gamma_n(t)) \, dt  \ = \ -4 \arctan \left( \frac{x}{\sqrt{\cos(x)}}
\right) \begin{array}{|l} x=\pi/2+2 \pi n \\ \\ x=-\pi/2 + 2 \pi n
\end{array}  \] which is zero when $n \neq 0$ and it is $-
4 \pi$ for $\gamma_0$. Each one of the other ovals of $f(x,y)=0$,
$\gamma_n$ with $n \neq 0$, belongs to the period annulus of a
center as it can be shown from the fact that the function
$H(x,y)=f(x,y) (x^2+y^2) \exp \{2 \arctan (y/x) \}$ is a first
integral for system (\ref{nalc}). Our result can be applied for
any of the periodic orbit $\gamma_n$ of this example.
\par
When considering a polynomial system, as far as the authors know,
only algebraic limit cycles are known in this implicit way. A
limit cycle is said to be {\em algebraic} if its points belong to
an invariant algebraic curve.
\newline

The paper is organized as follows. Section \ref{sect2} contains
the statement and proof of the main result of this work, i.e.
Theorem \ref{threl}. Using Theorem \ref{threl}, in Section
\ref{sect3}, we will show that all the known algebraic limit
cycles of a quadratic system are hyperbolic.

\section{Main result \label{sect2}}

\begin{theorem}
Let us consider a system {\rm (\ref{eq1})} and $\gamma(t)$ a
periodic orbit of period $T>0$. Assume that $f:\mathcal{U}
\subseteq \mathbb{R}^2 \to \mathbb{R}$ is an invariant curve with
$\gamma \subseteq \{ (x,y) \, | \, f(x,y) = 0\}$ and let $k(x,y)$
be the $\mathcal{C}^1$ function given in {\rm (\ref{ci})}. We
assume that if $p \in \mathcal{U}$ is such that $f(p)=0$ and
$\nabla f(p)=0$, then $p$ is a singular point of system {\rm
(\ref{eq1})}. Then,
\begin{equation}
\int_{0}^T  k (\gamma(t)) \, dt \  = \  \int_{0}^T
\div(\gamma(t)) \, dt. \label{eqr}
\end{equation}
\label{threl}
\end{theorem}
In order to prove Theorem \ref{threl}, we need to recall the
definition and some properties of the Poincar\'e map. Let us
consider $\gamma$ a periodic orbit with minimal period $T>0$ for
system (\ref{eq1}) and $p_0 \in \gamma$. Let $\mathcal{U}_\gamma
\subseteq \mathcal{U}$ be a neighborhood of $\gamma$ not
containing any singular point and $\Sigma = \{ q \in
\mathcal{U}_\gamma  \ | \ (q-p_0) \cdot \mathbf{F}(p_0) = 0 \}$,
where $\cdot$ denotes the scalar product between the vectors
$q-p_0$ and $\mathbf{F}(p_0)$.
\par As stated and proved in pages 210 and 211 in \cite{Perko}, we
have that there exists a $\delta >0$ and a unique function $\tau:
\Sigma \to \mathbb{R}$, which is defined continuously and
differentiable for any $q \in \Sigma \cap \mathcal{B}_\delta
(p_0)$ such that $\tau(p_0)=T$ and $\Phi_{\tau(q)}(q) \in \Sigma$.
As usual, $\mathcal{B}_\delta (p_0)$ is the ball of center $p_0$
and radius $\delta$. Then, for any $q \in \Sigma \cap
\mathcal{B}_\delta (p_0)$, the function $\mathcal{P}(q) =
\Phi_{\tau(q)}(q)$ is called the {\em Poincar\'e map} for $\gamma$
at $p_0$. It is clear that fixed points of the Poincar\'e map,
$\mathcal{P}(q)=q$, give rise to periodic orbits for system
(\ref{eq1}). Moreover, it can be shown that $\mathcal{P}: \Sigma
\to \Sigma$ is a $\mathcal{C}^1$ diffeomorphism.
\par
We consider $\Sigma$ as a subset of $\mathcal{U} \subseteq
\mathbb{R}^2$, so $\mathcal{P}$ is considered as a planar function
from $\Sigma \subset \mathbb{R}^2$ to $\mathbb{R}^2$. Hence, we
notice that the derivative of $\mathcal{P}$ at $p_0$, which is a
point in $\Sigma$, can be represented by a $2 \times 2$ matrix,
which we denote by $D \mathcal{P} (p_0)$. The following theorem,
stated and proved in \cite{Andronov2} page 118, is very useful to
establish the stability of $\gamma$.
\begin{theorem} Let $\mathbf{v}$ be a non-null vector normal to $\mathbf{F}(p_0)$. Then,
\begin{equation}
\mathbf{v} \cdot D\mathcal{P} (p_0) = \exp \left (\int_{0}^T
\div(\gamma(t)) \, dt \right) \mathbf{v}. \label{eqvn}
\end{equation}
\label{thDP}
\end{theorem}
Theorem \ref{thDP} is proved by using the variational equations of
first order related to system (\ref{eq1}). If $\Phi_t(x,y)$ is the
flow related to the vector field $\mathbf{F}(x,y)$, we have that
$\frac{d}{dt} \left(D\Phi_t(x,y)\right) = D\mathbf{F}(\Phi_t(x,y))
\cdot D\Phi_t(x,y)$ with the initial condition $D
\Phi_t(x,y)_{|t=0} =I$, where $D$ means the differential with
respect to the point $(x,y)$ and $I$ is the identity matrix. These
equations with respect to the matrix $D\Phi_t(x,y)$ are the
variational equations of first order. Since
$\mathcal{P}(q)=\Phi_{\tau(q)}(q)$, the solution of the
variational equations of first order allows the computation of
$D\mathcal{P} (p_0)$ in a point $p_0 \in \gamma$.
\par
In order to show that the stability of $\gamma$ is determined by
the value of $\mathbf{v} \cdot D\mathcal{P} (p_0)$, as stated in
Theorem \ref{thdiv}, we consider the {\em displacement function}
and we follow the reasoning of page 213 in \cite{Perko}. For any
$q \in \Sigma \cap \mathcal{B}_\delta (p_0)$, we have that $q =
p_0 + s \mathbf{v}$, with $s \in (-\delta/|\mathbf{v}| ,
\delta/|\mathbf{v}|)$. Since $\mathcal{P}(q) \in \Sigma$, we have
that given $s \in (-\delta/|\mathbf{v}|, \delta/|\mathbf{v}|)$,
there exists a $\sigma(s) \in \mathbb{R}$ such that
$\mathcal{P}(p_0 + s \mathbf{v}) = p_0 + \sigma(s) \mathbf{v}$.
Therefore, we have defined a $\mathcal{C}^1$ function $\sigma:
(-\delta/|\mathbf{v}|, \delta/|\mathbf{v}|) \to \mathbb{R}$ and
the displacement function is given by $\mathrm{d}:
(-\delta/|\mathbf{v}|, \delta/|\mathbf{v}|) \to \mathbb{R}$ with
$\mathrm{d}(s)= \sigma(s) - s$. It is clear that
$\mathrm{d}(0)=0$, $\mathrm{d}'(s)= \sigma'(s) -1$ and $\mathbf{v}
\cdot D\mathcal{P} (p_0 + s \mathbf{v})= \sigma'(s) \mathbf{v}$.
Since $\mathrm{d}(s)$ is $\mathcal{C}^1$, we have that the sign of
$\mathrm{d}'(s)$ coincides with the sign of $\mathrm{d}'(0)$ for
$|s|$ sufficiently small as long as $\mathrm{d}'(0) \neq 0$. By
the mean value theorem, we have that given $|s|$ sufficiently
small there exists a $\xi \in (0, s)$ such that $\mathrm{d}(s)=
\mathrm{d}'(\xi) s$. Therefore, if $\mathrm{d}'(0)
> 0$, we have that $\mathrm{d}(s)>0$ for $s>0$ and $\mathrm{d}(s)
< 0$ for $s<0$, which implies that the periodic orbit $\gamma$ is
an unstable limit cycle. Similar reasonings show that if
$\sigma'(0) >1$ then $\gamma$ is an unstable limit cycle and if
$\sigma'(0) <1$ then $\gamma$ is a stable limit cycle. Theorem
\ref{thdiv} clearly follows from Theorem \ref{thDP} and the fact
that $\sigma'(0) \mathbf{v} = \mathbf{v} \cdot D\mathcal{P}
(p_0)$.

\begin{lemma}
Let us consider a system {\rm (\ref{eq1})} and let $f:\mathcal{U}
\subseteq \mathbb{R}^2 \to \mathbb{R}$ be a non-null
$\mathcal{C}^1(\mathcal{U})$-function. There exists a
$\mathcal{C}^1$ function $k(x,y)$ such that $\nabla f (q) \cdot
\mathbf{F}(q)=k(q) f(q)$ for any $q \in \mathcal{U}$ if, and only
if, for any $q \in \mathcal{U}$ and any $t \in \mathbb{R}$ such
that $\Phi_t(q) \in \mathcal{U}$, the following identity is
satisfied:
\begin{equation}
f(\Phi_t(q)) = f(q) \ \exp \left( \int_{0}^t k (\Phi_s(q)) \, d s
\right). \label{eq5}
\end{equation}
\label{lemic}
\end{lemma}
\par {\em Proof.} Assume that $\nabla f (q) \cdot \mathbf{F}(q)=k(q)
f(q)$ for any $q \in \mathcal{U}$. We fix a point $q \in
\mathcal{U}$ and we define $\varphi(t)=f(\Phi_t(q))$ for any $t
\in \mathbb{R}$ such that $\Phi_t(q) \in \mathcal{U}$. We have
that $t$ belongs to an open interval $(-\epsilon_q, \epsilon_q)$
with $\epsilon_q >0$ (and it may be that $\epsilon_q = + \infty$).
We have, using some of the properties of the flow and the fact
$\dot{f}(\Phi_t(q)) = k(\Phi_t(q)) \, f(\Phi_t(q))$, that:
\[
\dot{\varphi}(t) =  \nabla f (\Phi_t(q)) \cdot \frac{d \Phi_t}{dt}
(q)  =  \nabla f (\Phi_t(q)) \cdot \mathbf{F} (\Phi_t(q))
 =  \dot{f} (\Phi_t(q))  =  k (\Phi_t(q)) \
\varphi(t).
\]
We deduce that $\frac{d \varphi}{dt}(t) = k (\Phi_t(q)) \
\varphi(t)$ and $\varphi(0)=f(q)$. Solving this linear equation in
the function $\varphi(t)$ we get $\varphi(t) = f(q) \exp \left(
\int_{0}^t k (\Phi_s(q)) \, d s \right)$. As we can consider the
same reasoning for any $q \in \mathcal{U}$, we obtain identity
(\ref{eq5}). The reciprocal is proved by the same reasoning. \bbox
\par
\begin{lemma}
Let us consider a system {\rm (\ref{eq1})} and $\gamma(t)$ a
periodic orbit of period $T>0$. Assume that $f:\mathcal{U}
\subseteq \mathbb{R}^2 \to \mathbb{R}$ is an invariant curve with
$\gamma \subseteq \{ (x,y) \, | \, f(x,y) = 0\}$ and let $k(x,y)$
be the $\mathcal{C}^1$ function given in {\rm (\ref{ci})}. Take
any $p_0$ in $\gamma$. Then,
\begin{equation}
\nabla f(p_0) \cdot D\mathcal{P} (p_0) = \exp \left( \int_{0}^T k
(\gamma(t)) \, dt \right) \nabla f (p_0). \label{eq4}
\end{equation}
\label{lemrel}
\end{lemma}
\par {\em Proof.} We consider the Poincar\'e map defined in an interval
of the straight line $\Sigma$ containing $p_0$, $\mathcal{P}(q)=
\Phi_{\tau(q)}(q)$. Since $f(x,y)=0$ is an invariant curve defined
in $\mathcal{U} \subseteq \mathbb{R}^2$, it is clear that for any
$q \in \mathcal{U}$ and any $t \in \mathbb{R}$ such that
$\Phi_t(q) \in \mathcal{U}$, identity (\ref{eq5}) is satisfied as
proved in Lemma \ref{lemic}. Hence,
\[ f(\mathcal{P}(q)) = f(q) \exp \left( \int_{0}^{\tau (q)}
k(\Phi_s(q)) \, d s \right), \] and differentiating this identity
with respect to $q$ we get
\[ \begin{array}{l}
\nabla f (\mathcal{P}(q)) \cdot D\mathcal{P} (q)  =   \exp \left(
\displaystyle \int_{0}^{\tau(q)} \displaystyle k(\Phi_s(q)) \, ds
\right) \nabla f (q) \, +
\\ \\ f(q)  \exp \left( \displaystyle \int_{0}^{\tau(q)} \displaystyle
k(\Phi_s(q)) ds \right) \! \! \left[ \displaystyle
\int_{0}^{\tau(q)} \left( \displaystyle
 \nabla k \right)
(\Phi_s(q)) \cdot D \Phi_s(q) \, ds + \displaystyle
k(\mathcal{P}(q)) \nabla \tau(q) \right],
\end{array} \]
where $D \mathcal{P}(q)$ and $D \Phi_s(q)$ stand for the Jacobian
matrix with respect to $q$ of the functions $\mathcal{P}$ and
$\Phi_s$, respectively, in the point $q$.
\par
We evaluate the previous identity in $q=p_0$, taking into account
that $f(p_0)=0$ and $\tau(p_0)=T$, and we get identity
(\ref{eq4}). \bbox
\newline

{\em Proof of Theorem {\rm \ref{threl}}.} The vector $\nabla f
(p_0)$ is a non-null vector that is normal to the vector
$\mathbf{F}(p_0)$ since $f(x,y)=0$ is an invariant curve that
contains $\gamma$, and $p_0 \in \gamma$. The fact of $\nabla f
(p_0)$ to be a non-null vector is ensured by the assumption that
if $p \in \mathcal{U}$ is such that $f(p)=0$ and $\nabla f(p)=0$,
then $p$ is a singular point of system {\rm (\ref{eq1})}. Since
$p_0$ belongs to a periodic orbit, it cannot be a singular point.
\par
Therefore, the vector $\mathbf{v}$ in the identity (\ref{eqvn}) of
Theorem \ref{thDP} can be replaced by $\nabla f(p_0)$. Using the
identity (\ref{eq4}) of Lemma \ref{lemrel}, we deduce that \[ \exp
\left ( \int_{0}^T \div (\gamma(t)) \, dt \right) \ = \ \exp \left
( \int_{0}^T k (\gamma(t)) \, dt \right),
\] from which (\ref{eqr}) follows. \bbox

\section{Hyperbolicity of the known algebraic limit cycles of quadratic systems \label{sect3}}

We consider the families of quadratic systems with algebraic limit
cycles known at the time of composition of this paper. These
families sweep all the algebraic limit cycles defined by
polynomials of degrees $2$ and $4$ for a quadratic system, as it
is proved in \cite{ChLlSo}. In \cite{Evdokimenko1, Evdokimenko2,
Evdokimenko3}, it is shown that there are no algebraic limit
cycles of degree $3$ for a quadratic system. See \cite{ChLlMo} for
a short proof. In \cite{ChLlSw}, two examples of quadratic systems
with an algebraic limit cycle of degree $5$ and $6$ are described.
We study the hyperbolicity of all these limit cycles.
\newline

The following result is due to Ch'in Yuan-sh\"{u}n \cite{Chin} and
characterizes the algebraic limit cycles of degree $2$ for a
quadratic system.
\begin{theorem} {\sc \cite{Chin}} \ If a quadratic system has an
algebraic limit cycle of degree $2$, then after an affine change
of variables, the limit cycle becomes the circle
\begin{equation}
\Gamma:= x^2+y^2-1=0. \label{f00} \end{equation} Moreover,
$\Gamma$ is the unique limit cycle of the quadratic system which
can be written in the form
\begin{equation}
\begin{array}{lll}
\dot{x} & = & -y\, (ax \, +\, by\, +\, c)-(x^2+y^2-1), \\
\dot{y} & = & x\, (ax\, +\, by\, +\, c),
\end{array}
\label{cicle2}
\end{equation}
with $a \neq 0$, $c^2 + 4 (b+1) >0$ and $c^2 > a^2 + b^2$.
\label{thchin}
\end{theorem}
We summarize the four families of algebraic limit cycles of degree
$4$ for quadratic systems in the following result, which is stated
and proved in \cite{ChLlSo}. We remark that these families were
encountered previously to the work \cite{ChLlSo}, but in this work
it was shown that there are no other algebraic limit cycle of
degree $4$ for a quadratic system. System (\ref{cicle4yab}) was
first described in \cite{Yablonskii}, system (\ref{cicle4fil}) in
\cite{Filipstov}, system (\ref{cicle4ch1}) in \cite{Chavarriga}
and system (\ref{cicle4ch2}) in \cite{ChLlSo}.
\begin{theorem} {\sc \cite{ChLlSo}} \
After an affine change of variables the only quadratic systems
having an algebraic limit cycle of degree $4$ are
\begin{itemize}
\item[{\rm (a)}] Yablonskii's system
\begin{equation}
\begin{array}{lll}
\dot{x} & = & -4  a  b  c  x - (a \, + \, b)  y + 3(a \, +  \, b) c  x^2 \, + \, 4  x  y, \\
\dot{y} & = & (a \, + \, b)  a  b  x - 4  a  b  c  y + (4  a  b
c^2 - \frac{3}{2} \, (a \, + \, b)^2 + 4  a  b)   x^2
\\ & &
 + 8 (a \, + \, b)  c  x  y + 8  y^2,
\end{array}
\label{cicle4yab}
\end{equation}
with $a b c \neq 0$, $a \neq b$, $a b>0$ and $4 c^2 (a - b)^2 + (3
a - b)(a - 3 b) < 0$. \\ This system has the invariant algebraic
curve
\begin{equation}
(y \, + \, c \, x^2)^2 + \, x^2 \, (x \, - \, a) (x \, - \, b) =0,
\label{f01}
\end{equation}
whose oval is a limit cycle for system {\rm (\ref{cicle4yab})}.

\item[{\rm (b)}] Filipstov's system
\begin{equation}
\begin{array}{lll}
\dot{x} & = & 6 \, (1 \, + \, a) \, x \, + \, 2 \, y \, - \, 6 \, (2 \, + \, a) \, x^2 \, + \, 12 \, x \, y, \\
\dot{y} & = & 15 \, (1 \, + \, a) \, y \, + \, 3 \, a \, (1 \, +
\, a) \, x^2 \, - \, 2 \, (9 \, + \, 5 \, a) \, x \, y \, + \, 16
\, y^2,
\end{array}
\label{cicle4fil}
\end{equation}
with $0< a < \frac{3}{13}$. This system has the invariant
algebraic curve
\begin{equation}
3(1 \, + \, a)(a \, x^2 \, + \, y)^2 \, + \, 2 \, y^2(2 \, y \, -
\, 3(1 \, + \, a)x)=0, \label{f02}
\end{equation}
whose oval is a limit cycle for system {\rm (\ref{cicle4fil})}.

\item[{\rm (c)}] Chavarriga's system
\begin{equation}
\begin{array}{lll}
\dot{x} & = & 5 \, x \, + \, 6 \, x^2 \, + \, 4(1 \, + \, a) \, x \, y \, + \, a \, y^2, \\
\dot{y} & = & x \, + \, 2 \, y \, + \, 4 \, x \, y \, + \, (2 \, +
\, 3 \, a) \, y^2,
\end{array}
\label{cicle4ch1}
\end{equation}
with $\frac{-71 + 17 \sqrt{17}}{32} < a <0$ has the invariant
algebraic curve
\begin{equation}
x^2 \, + \, x^3 \, + \, x^2 \, y \, + \, 2 \, a \, x \, y^2 \, +
\, 2 \, a \, x \, y^3 \, + \, a^2 \, y^4 =0, \label{f03}
\end{equation}
whose oval is a limit cycle for system {\rm (\ref{cicle4ch1})}.

\item[{\rm (d)}] Chavarriga, Llibre and Sorolla's system
\begin{equation}
\begin{array}{lll}
\dot{x} & = & 2 \, (1 \, + \, 2 \, x \, - \, 2 \, a \, x^2 \, + \, 6 \, x \, y), \\
\dot{y} & = & 8 \, - \, 3 \, a \, - \, 14 \, a \, x \, - \, 2 \, a
\, x \, y \, - \, 8 \, y^2,
\end{array}
\label{cicle4ch2}
\end{equation}
with $0 < a < \frac{1}{4}$ has the invariant algebraic curve
\begin{equation}
\frac{1}{4} \, + \, x \, - \, x^2 \, + \, a \, x^3 \, + \, x \, y
\, + \, x^2 \, y^2 =0, \label{f04}
\end{equation}
whose oval is a limit cycle for system {\rm (\ref{cicle4ch2})}.
\end{itemize}
\label{thcl4}
\end{theorem}

In a work due to C. Christopher, J. Llibre and G. \'{S}wirszcz
\cite{ChLlSw} two families of quadratic systems with an algebraic
limit cycle of degrees five and six, respectively, are given.
These two families are constructed by means of a birrational
transformation of system (\ref{cicle4ch2}). A {\em birrational
transformation} is a rational change of variables such that its
inverse is also rational. Moreover, they prove that there is also
a birrational transformation which converts Yablonskii's system
(\ref{cicle4yab}) into the system with a limit cycle of degree
$2$, that is, system (\ref{cicle2}).
\par
The fact of the limit cycle of degree $2$ being hyperbolic is
stated in \cite{Ye} (see pages 256--258) following the proof of
\cite{Chin}. As a consequence, and taking into account the
forthcoming Lemma \ref{lemdivg}, one of the limit cycles of degree
$4$ (the one due to Yablonskii) is also hyperbolic, because this
limit cycle of degree $4$ is birrationally equivalent to the one
of degree $2$, as it is shown in \cite{ChLlSw}. Our contribution
is the proof of the hyperbolicity of the other known limit cycles
of quadratic systems.

\begin{lemma}
Let us consider a differential system {\rm (\ref{eq1})} and a
change of variables $x=F(u,v)$ and $y=G(u,v)$, where $F,G$ are
$\mathcal{C}^2$ functions in $\mathcal{U}$. We denote by $
\dot{u}= R(u,v)$, $\dot{v}=S(u,v)$ the transformed differential
system. Let
\[ J(u,v):= \frac{\partial F}{\partial u} (u,v) \, \frac{\partial
G}{\partial v} (u,v) \, - \, \frac{\partial F}{\partial v} (u,v)
\, \frac{\partial G}{\partial u} (u,v), \] be the jacobian of the
transformation. Then,
\begin{equation}
\begin{array}{l}
\displaystyle \frac{\partial P}{\partial x} (F(u,v),G(u,v)) \, +
\,  \displaystyle \frac{\partial Q}{\partial y} (F(u,v),G(u,v)) \,
= \ \displaystyle \frac{\partial R}{\partial u} (u,v) \, +
 \displaystyle \frac{\partial S}{\partial v} (u,v) \, + \vspace{0.2cm} \\ \qquad \qquad + \,
 \displaystyle \frac{1}{J(u,v)} \left(  \displaystyle  \frac{\partial J}{\partial u} (u,v) \,
R(u,v) \, + \,  \displaystyle \frac{\partial J}{\partial v} (u,v)
\, S(u,v) \right).
\end{array} \label{eqdiv}
\end{equation} \label{lemdiv}
\end{lemma}
\par
\vspace{-0.2cm} Lemma \ref{lemdiv} is a computational result whose
proof is clear after some easy manipulations. We use it to prove
the following result which states that the value of the integral
of the divergence on the limit cycle does not change under
transformations of dependent variables.
\begin{lemma}
Let us consider a differential system (\ref{eq1}) with a periodic
orbit $\gamma$ of period $T>0$ and a change of variables
$x=F(u,v)$ and $y=G(u,v)$ which is well-defined in a neighborhood
of $\gamma$. We denote by $\dot{u}= R(u,v)$, $\dot{v}=S(u,v)$ the
transformed differential system and by $\vartheta$ the
corresponding periodic orbit. Then,
\[ \int_{0}^{T}  \left( \frac{\partial P}{\partial x} + \frac{\partial
Q}{\partial y}\right)(\gamma(t)) \, dt = \int_{0}^{T} \left(
\frac{\partial R}{\partial u} + \frac{\partial S}{\partial v}
\right)(\vartheta(t)) \, dt . \] \label{lemdivg}
\end{lemma}
\par
\vspace{-0.2cm} {\em Proof.} Using the same notation as in Lemma
\ref{lemdiv}, we have that the integral $ \int_{0}^{T}  \left(
\frac{\partial P}{\partial x} + \frac{\partial Q}{\partial y}
\right) (\gamma(t)) \, dt $ becomes, under the transformation of
dependent variables $x=F(u,v)$ and $y=G(u,v)$,
\[
\int_0^T \left( \frac{ \partial P}{\partial x} (F(u,v),G(u,v)) +
\frac{\partial Q}{\partial y} (F(u,v),G(u,v)) \right)
(\vartheta(t)) \, dt
\]
which, by Lemma \ref{lemdiv}, equals to
\[
\begin{array}{l} \displaystyle
\int_{0}^T \left( \frac{\partial R}{\partial u} (u,v) +
\frac{\partial S}{\partial v} (u,v) \right)(\vartheta(t))  \, dt
\, + \vspace{0.2cm}
\\ \displaystyle \
+ \int_{0}^T \frac{1}{J(u,v)} \left( \frac{\partial J}{\partial u}
(u,v) R(u,v) + \frac{\partial J}{\partial v} (u,v) S(u,v)
\right)(\vartheta(t)) \, dt. \end{array}
\]
We notice that the integrand of the second integral in the former
expression can be rewritten as $d(J(u,v))/J(u,v)$ and, since the
change of variables is well defined in a neighborhood of $\gamma$,
we have that this expression is a well defined, exact $1$-form
which is integrated over the closed curve $\vartheta$, so
$\oint_{\vartheta} d(J(u,v))/J(u,v) =0$. \bbox
\newline

Therefore, in order to prove that all these families of limit
cycles are hyperbolic, we only need to study the stability of the
limit cycles of systems (\ref{cicle4fil}), (\ref{cicle4ch1}) and
(\ref{cicle4ch2}). The hyperbolicity of the two limit cycles
described in \cite{ChLlSw} is shown by the fact that they are
birrationally equivalent to (\ref{cicle4ch2}).

\begin{theorem}
Each one of the limit cycles of systems {\rm (\ref{cicle4fil}),
(\ref{cicle4ch1})} and {\rm (\ref{cicle4ch2})} is hyperbolic.
\label{thhyp}
\end{theorem}
\par
{\em Proof.} \noindent In order to prove the hyperbolicity of the
limit cycles of systems (\ref{cicle4fil}), (\ref{cicle4ch1}) and
(\ref{cicle4ch2}) we use the same process for all of them. These
systems depend on a parameter $a$ which belongs to a certain open
interval when the limit cycle $\gamma$ exists. We denote by $T>0$
the period of the limit cycle and by $\mathcal{D}(a)$ the value of
the integral $\int_0^T \div(\gamma(t))\, dt$ for any value of the
parameter for which the limit cycle exists. This value decides the
hyperbolicity character of the limit cycle $\gamma$ in the system
with parameter $a$. By virtue of Lemma \ref{lemdivg}, we may
consider any birrational transformation of these systems well
defined in a neighborhood of the limit cycle and we may consider
the transformed system instead of the previous one because the
value of the integral $\int_0^T \div(\gamma(t))\, dt$ does not
change.
\par
Using Theorem \ref{threl}, we have that:
\[ \mathcal{D}(a) = \int_0^T \div(\gamma(t)) \, dt + w \left( \int_0^T
\div(\gamma(t)) \, dt - \int_0^T k(\gamma(t)) \, dt \right), \]
where $k$ is the cofactor of the invariant algebraic curve
containing the limit cycle and $w$ is any real number.
\par
We show that the function $\mathcal{D}(a)$ has no zero when $a$
belongs to the interval of existence of limit cycle by choosing an
adequate $w \in \mathbb{R}$ and parameterizing the limit cycle
$\gamma$. The way of choosing the adequate value of $w$ is purely
heuristic, although we expect that this choice is related to some
geometric property. We find it very surprising that it is possible
to choose $w=-3$ for each one of the three families of systems.
\newline

\noindent {\em Hyperbolicity of the limit cycle given by the
algebraic curve {\rm (\ref{f04})} for system {\rm
(\ref{cicle4ch2}).}}
\newline

The stability of the limit cycle $\gamma$ is given by the
following function of the parameter $a$ of the system,
$\mathcal{D}(a) := \int_{0}^{T} \div(\gamma(t)) \, dt$, where
$\div(x,y)=2(2-5 a x - 2 y)$ is the divergence of system
(\ref{cicle4ch2}) and $T>0$ the period of the limit cycle. Theorem
\ref{threl} gives:
\[ \int_{0}^T \div(\gamma(t)) \, dt = \int_{0}^T k(\gamma(t)) \, dt,
\] where $k(x,y)= 4(2-3ax+2y)$ is the cofactor of the invariant
algebraic curve (\ref{f04}). So, given any real number $w$, we
have that:
\begin{eqnarray*}
\mathcal{D}(a) & = & \int_{0}^T \div(\gamma(t)) \, dt + w
\int_{0}^T (\div-k)(\gamma(t)) \, dt \\ & = & \int_{0}^T \left(
(1+w) \div - w k \right) (\gamma(t)) \, dt .
\end{eqnarray*}
We consider the following parameterization of the oval of the
algebraic curve (\ref{f04}):
\begin{equation}
x(\tau)=\tau, \quad y_{\pm}(\tau) = \frac{ - 1 \pm 2 \sqrt{(-a)
\tau (\tau-\tau_1) (\tau- \tau_2)}}{2 \tau}, \label{Dparch2}
\end{equation}
where $\tau_1= \frac{ 1-\sqrt{1-4a}}{2a}$, $\tau_2= \frac{1+
\sqrt{1-4a}}{2a}$ and the parameter $\tau \in (\tau_1, \tau_2)$.
The positive sign $y_{+}(\tau)$ gives a half of the oval and the
negative sign $y_{-}(\tau)$ the other half. One of the endpoints
of both parameterizations is $(x_1,y_1)=(\frac{
1-\sqrt{1-4a}}{2a}, -\frac{1+\sqrt{1-4a}}{4})$ and the other
endpoint is $(x_2,y_2)=(\frac{ 1+\sqrt{1-4a}}{2a},
\frac{-1+\sqrt{1-4a}}{4})$. We have that the vector field in
$(x_1,y_1)$ is $(0, 6 \sqrt{1-4a})$ and in $(x_2,y_2)$ is $(0, - 6
\sqrt{1-4a})$, so the flow on the limit cycle is clockwise. The
line $2 a x =1$ cuts the limit cycle in two points with ordinates
$\pm \sqrt{\frac{1-4a}{2}} -a$, which are given respectively by
$y_{\pm}(1/2a)$. We have the following relation between the
differentials: $d \tau = P (x (\tau), y_{\pm}(\tau)) \, dt$ where
$P(x,y)= 2  (1  +  2 x - 2  a x^2 + 6  x  y)$. Then,
\begin{eqnarray*}
\mathcal{D}(a)  & = & \displaystyle \int_{0}^T \displaystyle
\left( (1+w) \, \div - w k \right) (\gamma(t)) \, dt \vspace{0.2cm} \\
\vspace{0.2cm}  & = & \displaystyle \int_{\tau_1}^{\tau_2} \! \!
\! \displaystyle \left( \displaystyle \frac{\left( (1+w) \, \div -
w k \right)}{P} \right) (\tau, y_{+} (\tau)) \, d \tau +
\\ \vspace{0.2cm} & & \qquad  + \displaystyle \int_{\tau_2}^{\tau_1} \! \! \! \displaystyle \left(
\displaystyle \frac{\left( (1+w)\,  \div - w k
\right)}{P} \right) (\tau, y_{-} (\tau)) \, d \tau \\
\vspace{0.2cm} & = & \displaystyle \int_{\tau_1}^{\tau_2} \left[
\displaystyle  \left( \displaystyle \frac{\left( (1+w) \, \div - w
k
\right)}{P} \right) (\tau, y_{+} (\tau)) \, \right. \\
\vspace{0.2cm} & & \qquad \left. - \, \displaystyle \left(
\displaystyle \frac{\left( (1+w) \, \div - w k \right)}{P} \right)
(\tau, y_{-} (\tau)) \right] \, d \tau .
\end{eqnarray*}
For $w=-3$ and substituting by the parameterization, we get,
\[ \mathcal{D}(a) = 8 \int_{\tau_1}^{\tau_2} \frac{\sqrt{a \tau
(\tau- \tau_1) (\tau_2 - \tau)}}{\tau (1 + 8 \tau + a \tau^2)} \,
d \tau. \] Since $\tau_1>0$ and $\tau_2>\tau_1$ for any $a \in (0,
1/4)$ and the integrand $\frac{\sqrt{a \tau (\tau- \tau_1) (\tau_2
- \tau)}}{\tau (1 + 8 \tau + a \tau^2)}$ is also strictly positive
and well defined for any $\tau \in (\tau_1, \tau_2)$ and $a \in
(0, 1/4)$, we have that $\mathcal{D}(a)>0$ for all $a \in (0,
1/4)$, which implies that the limit cycle in system
(\ref{cicle4ch2}) is hyperbolic (and unstable).
\newline

\noindent {\em Hyperbolicity of the limit cycle given by the
algebraic curve {\rm (\ref{f02})} for system {\rm
(\ref{cicle4fil}).}} \newline

In order to simplify our computations, we consider the following
birrational change of the parameter, $a=3c/(4+5c)$, with inverse
$c=4a/(3-5a)$ and we have that $c \in (0,1/2)$.
\par
We consider the birrational change of variables $(x,y) \to (u,v)$
given by $(x,y)=(X(u,v),Y(u,v))$, where
\begin{eqnarray*}
X(u,v) & = & - \frac{2(1+2c)}{c^2 (1+u)^3} \, (c - 2(1+c)u + c u^2
- 2 \sqrt{1+2c}\, v), \vspace{0.2cm} \\
Y(u,v) & = & - \frac{12 (1+ 2 c)^2}{c^2 (4 + 5 c) (1 + u)^4} \, (c
- 2(1+c)u + c u^2 - 2 \sqrt{1+2c}\, v).
\end{eqnarray*}
The inverse of this change is given by $(u,v)=(U(x,y),V(x,y))$
with
\begin{eqnarray*}
U & = & \frac{6(1+2c) x}{(4+5c)y} -1 , \vspace{0.2cm} \\ V & = &
\frac{\sqrt{1+2c}}{(4+5c)^3 y^3} \bigg[ (1+2c) \big( 54 c^2 x^4 +
18 c (4 + 5c) x^2 y - 6 (4 + 5 c)^2 x y^2 \big) \bigg] +
\sqrt{1+2c}.
\end{eqnarray*}
The jacobian of this change of variables is: \[ \frac{\partial
U}{\partial x} \frac{\partial V}{\partial y} - \frac{\partial
U}{\partial y} \frac{\partial V}{\partial x} \ = \  \frac{324 \,
c^2 \, (1 + 2 c)^{5/2} \ x^4}{(4+5c)^4 \, y^5},
\] which can be seen to be well defined and different from zero in
all the points of the oval of the curve given in (\ref{f02}).
\par
We get a transformed system in which we reparameterize its time
$t$ multiplying by $c^2(4+5c)(1+u)^3/(12 (1+2c))$. This
reparameterization does not affect the direction of the flow on
the limit cycle. The new system reads for:
\begin{equation}
\begin{array}{lll}
\displaystyle \dot{u} & = & -2 u \, (cu + 4 +9c)
(cu^2-u+c)  - \displaystyle  2  \sqrt{1+2c}\, (cu^2- (4+5c)u + 2 c) v, \vspace{0.2cm} \\
\displaystyle \dot{v} & = & -c^2 \, \sqrt{1+2c} \, (u+1)^2(u-1) (3u+2)  \vspace{0.2cm}\\
& &  \displaystyle - (cu+4+9c) (3cu^2-2u+c)v \displaystyle + 2 \,
\sqrt{1+2c} \, (4+ 5c- 3cu) v^2,
\end{array}
\label{eqfil2}
\end{equation}
and the limit cycle is transformed to the real oval of the curve $
v^2 + u (c u^2 - u + c) =0$. This algebraic curve is invariant for
system (\ref{eqfil2}) with cofactor $k(u,v)=4  \sqrt{1+2c} \,
(4+5c-3cu)v - 2  (cu+9c+4)(3cu^2-2u+c)$. The divergence of system
(\ref{eqfil2}) is $\div(u,v)=2  \sqrt{1+2c} \, (12+15c-8cu) v - [
11 c^2 u^3+ c (28+81c) u^2 + (5 c^2 - 54 c - 24) u + 3 c (4 + 9 c)
]$.
\par
We consider the following parameterization of the oval of the
algebraic curve $v^2 + u (c u^2 - u + c) =0$:
\begin{equation}
u(\tau)=\tau, \quad v_{\pm}(\tau) = \pm \sqrt{c \tau (\tau-\tau_1)
(\tau_2 - \tau)}, \label{Dparfil}
\end{equation}
where $\tau_1= \frac{ 1-\sqrt{1-4c^2}}{2c}$, $\tau_2= \frac{1+
\sqrt{1-4c^2}}{2c}$ and the parameter $\tau \in (\tau_1, \tau_2)$.
We notice that for $c \in (0,1/2)$, we have that $0 < \tau_1 < 1 <
\tau_2$. The endpoints of both parameterizations are $(\tau_1,0)$
and $(\tau_2,0)$, and the vector field at the point $(\tau_i,0)$
is $(0,c^2 \sqrt{1+2c}\, (\tau_i+1)^2(3\tau_i+2)(1-\tau_i))$, for
$i=1,2$. Therefore, the flow on the limit cycle is clockwise. We
follow analogous arguments to the previous example to deduce that:
\begin{eqnarray*}
\mathcal{D}(c)  & = & \displaystyle \int_{\tau_1}^{\tau_2} \left[
\displaystyle  \left( \displaystyle \frac{\left( (1+w) \, \div - w
k
\right)}{P } \right) (\tau, v_{+} (\tau)) \, \right. \\
\vspace{0.2cm} & & \qquad \left. - \, \displaystyle \left(
\displaystyle \frac{\left( (1+w) \, \div - w k \right)}{P} \right)
(\tau, v_{-} (\tau)) \right] \, d \tau ,
\end{eqnarray*}
where $P(u,v)$ is the polynomial which defines $\dot{u}=P(u,v)$.
For $w=-3$ and substituting by the parameterization, we get
\[ \mathcal{D}(c) = 8 \sqrt{1+2c} \int_{\tau_1}^{\tau_2} \frac{\sqrt{c \tau
(\tau- \tau_1) (\tau_2 - \tau)}}{(\tau+1) (c \tau^2 + (17c+8) \tau
+ 4 + 8c)} \, \, d \tau.
\] Since $0<\tau_1<1<\tau_2$ for any $c \in (0, 1/2)$
and the integrand is strictly positive and well defined for any
$\tau \in (\tau_1, \tau_2)$ and $c \in (0, 1/2)$, we have that
$\mathcal{D}(c)>0$ for all $c \in (0, 1/2)$, which implies that
the limit cycle in system (\ref{eqfil2}) given by the real oval of
$v^2+ u (c u^2 - u + c) =0$ is hyperbolic (and unstable). Hence,
using Lemma \ref{lemdivg}, we have that the limit cycle given by
the oval of the curve (\ref{f02}) in system (\ref{cicle4fil}) is
hyperbolic and unstable.
\newline

\noindent {\em Hyperbolicity of the limit cycle given by the
algebraic curve {\rm (\ref{f03})} for system {\rm
(\ref{cicle4ch1}).}} \newline

We consider the following birrational change of variables $(x,y)
\to (u,v)$ with $(x,y)= (X(u,v),Y(u,v))$ where
\[ (X(u,v),Y(u,v)) \, = \, \frac{(-2,-2u)}{v+1+u+2 a u^2}. \] The
inverse of this change is $(u,v)=(U(x,y),V(x,y))$ with
\[ U(x,y)=\frac{y}{x}, \quad V(x,y)=-2a \, \frac{y^2}{x^2} -
\frac{(y+2)}{x} -1. \] This change of variables is well defined in
a neighborhood of the real oval of the curve (\ref{f03}) and its
jacobian is \[ \frac{\partial U}{\partial x} \frac{\partial
V}{\partial y} - \frac{\partial U}{\partial y} \frac{\partial
V}{\partial x} \, = \, -\frac{2}{x^3}. \] The algebraic curve
(\ref{f03}) is transformed to $v^2 + 4 a u^2 (u-1)- (u+1)^2=0$. We
consider the transformed system in which we make a
reparameterization of its time $t$ which consists on multiplying
by $v+1+u+2 a u^2$. This reparameterization reverses the direction
of the flow on the transformed limit cycle and the new system is
written:
\begin{equation}
\begin{array}{lll}
\displaystyle \dot{u} & = & (u+1)^2 - 4 a u^2(u-1) + (1-3u) v,
\vspace{0.2cm} \\ \displaystyle \dot{v} & = & 2 (u+1) (3
+u+2au-au^2) + (1+4au+u-6au^2) v - 5 v^2.
\end{array}
\label{eqch12}
\end{equation}
The divergence of this system is $\div(u,v) = 3(1+u) + 12 a u - 18
a u^2 - 13 v$. The algebraic curve $v^2 + 4 a u^2 (u-1)-
(u+1)^2=0$ is invariant for system (\ref{eqch12}) with cofactor
$k(u,v) = 2 (1+u+ 4 a u - 6 a u^2 - 5 v)$. The real oval of this
curve is a hyperbolic limit cycle for system (\ref{eqch12}) if,
and only if, the real oval of the curve (\ref{f03}) is a
hyperbolic limit cycle for system (\ref{cicle4ch1}), by Lemma
\ref{lemdivg}. We are going to compute the function
$\mathcal{D}(a) = \int_{0}^T \div(\gamma(t)) \, dt$, where
$\gamma$ is the real oval of the curve $v^2 + 4 a u^2 (u-1)-
(u+1)^2=0$. To do so, we parameterize this oval by:
\[ u(\tau)=\tau, \qquad v_{\pm}(\tau) = \pm \sqrt{ 4 a \tau^2 (1-\tau) +
( \tau+1)^2}, \] where $\tau$ takes values between $\tau_1$ and
$\tau_2$. The values $\tau_1$ and $\tau_2$ are the two smallest
roots of the polynomial $g(a, \tau):=4 a \tau^2 (1-\tau) + (
\tau+1)^2$ in $\tau$. We consider $a$ in the interval between $(17
\sqrt{17}-71)/32$ and $0$, which are the values of the parameter
for which the limit cycle exists. Since the coefficient of the
highest order term of $g(a,\tau)$ is $-4a$ which is strictly
negative, $g(a, -(3+\sqrt{17})/2)=4 (29 + 7\sqrt{17}) \left[ a -
(17 \sqrt{17}-71)/32 \right] >0$ and $g(a,-1)=8a<0$, we deduce
that $\tau_1 < -(3+\sqrt{17})/2 < \tau_2 < -1$. We denote by
$P(u,v)$ the polynomial which defines $\dot{u}=P(u,v)$ in system
(\ref{eqch12}). We consider the point with coordinates: \[
(u_0,v_0) \, = \, \left( -\frac{3+\sqrt{17}}{2}\, , 2\,
\sqrt{29+7\sqrt{17}} \, \sqrt{a+\frac{71-17 \sqrt{17}}{32}}
\right) \] and we have that $v_0^2-g(a,u_0)=0$ and $P(u_0,v_0) =
v_0 [ v_0 + (11+3\sqrt{17})/2 ]$. We deduce that $P(u_0,v_0) >0$
and that the flow on the limit cycle is clockwise. Using analogous
arguments as in the previous examples we conclude that:
\begin{eqnarray*}
\mathcal{D}(a)  & = & \displaystyle \int_{\tau_1}^{\tau_2} \left[
\displaystyle  \left( \displaystyle \frac{\left( (1+w) \, \div - w
k
\right)}{P(u,v)} \right) (\tau, v_{+} (\tau)) \, \right. \\
\vspace{0.2cm} & & \qquad \left. - \, \displaystyle \left(
\displaystyle \frac{\left( (1+w) \, \div - w k \right)}{P(u,v)}
\right) (\tau, v_{-} (\tau)) \right] \, d \tau .
\end{eqnarray*}
For $w=-3$ and substituting by the parameterization, we get
\[ \mathcal{D}(a) = 2 \int_{\tau_1}^{\tau_2} \frac{\sqrt{g(a,\tau)}}{(\tau-1) \,\tau\, (a \tau+2)} \, \, d \tau.
\] We have that $\tau_1< -(3+\sqrt{17})/2 <\tau_2<-1$ and that $\tau-1<0$ and $\tau<0<-2/a$ for any $\tau \in
(\tau_1, \tau_2)$ and for any $a \in ( (17 \sqrt{17}-71)/32 , 0)$.
Hence, the integrand is strictly positive and well defined for any
$\tau \in (\tau_1, \tau_2)$ and $a \in ( (17 \sqrt{17}-71)/32 ,
0)$. We deduce that $\mathcal{D}(a)>0$ for all $a \in ( (17
\sqrt{17}-71)/32 , 0)$, which implies that the limit cycle in
system (\ref{eqch12}) given by the real oval of $v^2 - g(a,u)=0$
is hyperbolic (and unstable). Thus, using Lemma \ref{lemdivg} and
the sign in the change of time, we have that the limit cycle given
by the oval of the curve (\ref{f03}) in system (\ref{cicle4ch1})
is hyperbolic and stable. \bbox
\newline

\section*{Appendix}

The aim of this appendix is to present some relations among
elliptic integrals which the authors obtained by using the
identity given by Theorem \ref{threl} for systems
(\ref{cicle4fil}) and (\ref{cicle4ch2}).
\par
Before the presented proof of Theorem \ref{thhyp}, the authors got
its proof for systems (\ref{cicle4fil}) and (\ref{cicle4ch2}) by
computing the corresponding integrals which give place to elliptic
integrals. The identity given in Theorem \ref{threl} was used to
encounter a Fuchs equation for the function $\mathcal{D}(a)$.
After some thorough analysis of this Fuchs equation, we deduce the
non-vanishing of the function $\mathcal{D}(a)$ for any value of
the parameter in which the limit cycle exists. We are not going to
give this proof, but we think that the relations among elliptic
integrals obtained by the former reasoning are interesting by
themselves. Hence, we give the identities obtained which, as far
as we know, do not appear in any book of tables of integrals and
relations between classical functions. On the other hand, we also
give the obtention of the Fuchs equation for the function
$\mathcal{D}(a)$ in the case of system (\ref{cicle4ch2}).

\subsection*{Identities among elliptic integrals}

The functions involved in this subsection are the complete
elliptic integrals of first, second and third kinds, denoted by
$\mathrm{K}(\omega)$, $\mathrm{E}(\omega)$ and
$\mathrm{\Pi}(\kappa, \omega)$, respectively. We recall the
definition of these functions:
\begin{eqnarray*}
\mathrm{K}(\omega) & = & \displaystyle \int_{0}^{\pi/2}
\displaystyle \frac{ d \theta}{\sqrt{1-\omega \sin^2(\theta)}}  =
 \displaystyle \int_{0}^{1}  \displaystyle \frac{d
t}{\sqrt{(1-t^2) (1- \omega t^2)}}, \\ \mathrm{E}(\omega) & = &
\displaystyle \int_{0}^{\pi/2} \sqrt{1-\omega \sin^2(\theta)} \, d
\theta  =   \displaystyle \int_{0}^{1}
 \displaystyle \frac{\sqrt{1-\omega t^2}}{\sqrt{(1-t^2)}} \, d t,
 \\
\mathrm{\Pi} (\kappa, \omega) & = &  \displaystyle
\int_{0}^{\pi/2} \displaystyle  \frac{ d \theta}{(1-\kappa \sin^2
(\theta)) \sqrt{1-\omega \sin^2(\theta)}}  =   \displaystyle
\int_{0}^{1}  \displaystyle  \frac{d t}{(1- \kappa t^2)
\sqrt{(1-t^2) (1- \omega t^2)}} ,
\end{eqnarray*}
and their derivatives:
\begin{eqnarray*}
\mathrm{K}'(\omega) & = & \frac{1}{2(1-\omega)\omega} \,
\mathrm{E}(\omega) - \frac{1}{2 \omega} \, \mathrm{K} (\omega), \\
\mathrm{E}'(\omega) & = & \frac{1}{2 \omega} \, \mathrm{E}
(\omega) - \frac{1}{2 \omega} \mathrm{K} (\omega), \\
\frac{\partial \mathrm{\Pi} (\kappa, \omega)}{\partial \kappa}& =
& \frac{1}{2 \kappa (\kappa-1)} \, \mathrm{K}(\omega) + \frac{1}{2
(\kappa-1) (\omega - \kappa)} \, \mathrm{E}(\omega) +
\frac{\kappa^2- \omega}{2 \kappa (\kappa-1) (\omega - \kappa)} \,
\mathrm{\Pi} (\kappa, \omega), \\
\frac{\partial \mathrm{\Pi} (\kappa, \omega)}{\partial \omega} & =
& \frac{1}{2 (\kappa - \omega) (\omega -1)} \, \mathrm{E} (\omega)
+ \frac{1}{2 (\kappa - \omega)} \, \mathrm{\Pi} (\kappa, \omega) .
\end{eqnarray*}
\newline

We use the following parameterization of the oval of the algebraic
curve (\ref{f04}) to explicitly compute the integrals for the
system (\ref{cicle4ch2}). We parameterize the oval by:
\begin{equation}
x(\tau)=\tau, \quad y_{\pm}(\tau) = \frac{ - 1 \pm 2 \sqrt{(-a)
\tau (\tau-\tau_1) (\tau- \tau_2)}}{2 \tau}, \label{parch2}
\end{equation}
where $\tau_1= \frac{ 1-\sqrt{1-4a}}{2a}$, $\tau_2= \frac{1+
\sqrt{1-4a}}{2a}$ and the parameter $\tau \in (\tau_1, \tau_2)$.
The positive sign $y_{+}(\tau)$ gives a half of the oval and the
negative sign $y_{-}(\tau)$ the other half. The explicit
computation of the integrals for the system (\ref{cicle4ch2})
gives that the identity (\ref{eqr}) stated in Theorem \ref{threl}
reads for:
\begin{equation}
-9 \, \mathrm{K}(\omega_0) \, + \,  c_{+} \,
\mathrm{\Pi}(\omega_{+}, \omega_0) \, + \,  c_{-} \,
\mathrm{\Pi}(\omega_{-}, \omega_0) \equiv 0 , \label{relch2}
\end{equation}
which is valid for $a \in (0, 1/4)$, where \[ \omega_0 = \frac{2
\sqrt{1-4a}}{1+\sqrt{1-4a}}, \  \omega_{\pm} = \frac{2
\sqrt{1-4a}}{9 + \sqrt{1-4a} \pm 2 \sqrt{16-a}}, \  c_{\pm}=
\frac{9 - \sqrt{1-4a}}{2} \pm \sqrt{16-a} . \] The derivative of
the expression in (\ref{relch2}) with respect to $a$ gives place
to the same identity (\ref{relch2}). In fact, when computing the
derivative with respect to $a$ of the expression given in
(\ref{relch2}), using the described formulas of derivation for
these elliptic integrals, we get $-1/(1-4a + \sqrt{1-4a})$ times
the same expression (\ref{relch2}). This simple factor is
different from zero when $a \in (0,1/4)$.
\par
In the same way, we can explicitly compute the integrals involved
in the identity (\ref{eqr}) stated in Theorem \ref{threl} for the
system (\ref{cicle4fil}), via using the parameterization of the
oval of (\ref{f02}) given by:
\begin{equation}
\begin{array}{l}
x_{\pm}(\tau)  =  - \displaystyle \frac{2(1+2c)}{c^2(1+\tau)^3}
(c-2 \tau - 2 c \tau + c \tau^2 \pm 2
\sqrt{1+2c}\sqrt{\tau(-c+\tau-c \tau^2)}),
\\ \\ y_{\pm}(\tau)  =  \displaystyle \frac{-12(1+2c)^2}{c^2(4+5c)(1+ \tau)^4} (c
- 2\tau - 2 c \tau + c \tau^2 \pm 2
\sqrt{1+2c}\sqrt{\tau(-c+\tau-c \tau^2)}),
\end{array} \label{parfil} \end{equation}
where $\tau \in (\tau_1, \tau_2)$ with $\tau_1=\frac{1-
\sqrt{1-4c^2}}{2 c}$ and $\tau_2=\frac{1+ \sqrt{1-4c^2}}{2 c}$. It
is clear that $0< \tau_1 < 1 < \tau_2$ for $c \in (0,1/2)$. We
define $g(c,\tau)=-\tau(c-\tau+c \tau^2)=c \tau
(\tau-\tau_1)(\tau_2-\tau)$, which is strictly positive for all
$\tau \in (\tau_1,\tau_2)$. The explicit computation of the
integrals involved in the identity (\ref{eqr}) gives:
\begin{equation}
5 \, \mathrm{K}(\varsigma_0) \, + \,  C_{+} \,
\mathrm{\Pi}(\varsigma_{+}, \varsigma_0) \, + \,  C_{-} \,
\mathrm{\Pi}(\varsigma_{-}, \varsigma_0) \equiv 0 , \label{relfil}
\end{equation}
which is valid for $c \in (0, 1/2)$, where \[ \varsigma_0 =
\frac{2 \sqrt{1-4c^2}}{1+\sqrt{1-4c^2}}, \quad \varsigma_{\pm} =
\frac{2 \sqrt{1-4c^2}}{9 + 17 c + \sqrt{1-4c^2} \pm
\sqrt{64(1+2c)^2+c^2}}, \] \[ C_{\pm}= \frac{-2(24+47c \pm 3
\sqrt{64(1+2c)^2+c^2}}{9 + 17 c + \sqrt{1-4c^2} \pm
\sqrt{64(1+2c)^2+c^2}} . \] The derivative of the expression
(\ref{relfil}) with respect to $c$ gives place to the same
identity (\ref{relfil}).
\par
The authors have not been able to give an analogous identity
related to system (\ref{cicle4ch1}) due to the fact that the
corresponding integrals require much more computations to be
identified with the elliptic integrals.

\subsection*{Fuchs equation for $\mathcal{D}(a)$ in system
(\ref{cicle4ch2})}

In this part of the appendix we develop the way we obtained a
Fuchs equation for the function $\mathcal{D}(a)$ in system
(\ref{cicle4ch2}), via using the relation (\ref{relch2}). We think
that the fact of obtaining a Fuchs equation satisfied by this
function is interesting to further understand the stability of
algebraic limit cycles for polynomial systems. We obtained a
similar Fuchs equation for system (\ref{cicle4fil}), but we do not
state it because the equation itself does not give any further
information about the properties of system (\ref{cicle4fil}) and
the way it was obtained is completely analogous to the way
equation (\ref{fuchsch2}) for system (\ref{cicle4ch2}) is
obtained.
\par
Let us consider system (\ref{cicle4ch2}) and we parameterize the
oval which contains the limit cycle by (\ref{parch2}). Taking the
notation described in the previous subsection: \[ \omega_0 =
\frac{2 \sqrt{1-4a}}{1+\sqrt{1-4a}}, \ \omega_{\pm} = \frac{2
\sqrt{1-4a}}{9 + \sqrt{1-4a} \pm 2 \sqrt{16-a}}, \ c_{\pm}=
\frac{9 - \sqrt{1-4a}}{2} \pm \sqrt{16-a}, \] and \[ \mu=\sqrt{1+
\sqrt{1-4a}}, \quad b_{\pm}=2(4 \pm \sqrt{16-a})\, c_{\pm}, \] we
explicitly compute the value of $\mathcal{D}(a)$:
\begin{eqnarray*}
\mathcal{D}(a) & = & \int_{0}^{T} \div (\gamma(t)) \, dt  \\ & = &
\frac{\sqrt{2}}{\mu \sqrt{16-a}} \left[ -34 \sqrt{16-a} \,
\mathrm{K}(\omega_0) + b_{+} \, \mathrm{\Pi} (\omega_{+},
\omega_0) - b_{-} \, \mathrm{\Pi} (\omega_{-}, \omega_0) \right] .
\end{eqnarray*}
We compute the successive derivatives of $\mathcal{D}(a)$:
\begin{eqnarray*}
\mathcal{D}'(a) & = & \frac{-4 \sqrt{2}}{\mu \, (16-a)^{3/2} }
\bigg[ \sqrt{16-a} \, \mathrm{K}(\omega_0) \, + \, \frac{2 \mu^2
\, \sqrt{16-a}}{a} \, \mathrm{E}(\omega_0) \, +  \\ & &  - \,
c_{+} \, \mathrm{\Pi} (\omega_{+}, \omega_{0}) \, + \, c_{-} \,
\mathrm{\Pi} (\omega_{-}, \omega_{0}) \bigg], \\
\mathcal{D}''(a) & = & \frac{6 \sqrt{2}}{\mu \, (16-a)^{5/2}}
\Bigg[ \frac{(10a^2+33 a -64) \sqrt{16-a}}{3a (1-4a)} \,
\mathrm{K}(\omega_0) \, + \\ & & \frac{(73 a^2 - 420 a + 128)
\mu^2 \sqrt{16-a}}{6 a^2 (1-4a)} \, \mathrm{E}(\omega_0)  + \,
c_{+} \, \mathrm{\Pi}(\omega_{+}, \omega_{0}) \, - \, c_{-} \,
\mathrm{\Pi}(\omega_{-}, \omega_{0}) \Bigg], \\
\mathcal{D}'''(a) & = & \frac{- \sqrt{2}}{\mu (16-a)^{7/2}} \Bigg[
\frac{(180 a^4 + 1347 a^3 - 9685 a^2 + 25664 a - 4096)
\sqrt{16-a}}{a^2 (1-4a)^2} \, \mathrm{K} (\omega_0) \\ & & +
\frac{(1812 a^4 - 20259 a^3 + 102164 a^2 - 60544 a + 8192) \mu^2
\sqrt{16-a}}{2 a^3 (1 - 4 a)^2} \, \mathrm{E}(\omega_0) + \\ & & -
\, 15 c_{+} \, \mathrm{\Pi}(\omega_{+}, \omega_{0}) \, + \, 15
c_{-} \mathrm{\Pi}(\omega_{-}, \omega_{0}) \Bigg].
\end{eqnarray*}
By elimination of independent functions and using the identity
(\ref{relch2}) we obtain the following third order homogeneous
differential equation of Fuchs type for $\mathcal{D}(a)$:
\begin{equation}
\begin{array}{l}
8 (a-16) a (4 a -1) (17 a +8) \, \mathcal{D}'''(a) \, + \, 4 ( 612
a^3 - 4119 a^2 - 2600 a \\ + 512) \, \mathcal{D}''(a)
\vspace{0.2cm} + \, 6 (a-2) (289 a + 528) \, \mathcal{D}'(a) \, +
\, 3 (17 a + 64) \, \mathcal{D}(a)  = \, 0.
\end{array}
\label{fuchsch2}
\end{equation}
An easy computation shows that $\mathcal{D}(1/4) = 0$,
$\mathcal{D}'(1/4)= -8 \sqrt{2} \pi/9$ and $\mathcal{D}''(1/4)= 98
\sqrt{2} \pi/ 27$. Hence, equation (\ref{fuchsch2}) univocally
determines the function $\mathcal{D}(a)$ defined in $a \in (0,
1/4]$. A thorough analysis of the properties of $\mathcal{D}(a)$
gives that $\mathcal{D}(a) > 0$ for $a \in (0,1/4)$.
\par
We remark that using identity (\ref{relch2}) we get a Fuchs
equation of order $3$ for $\mathcal{D}(a)$. If we did not have
this relation, we would get an equation of order $4$, which would
make the analysis of properties much more difficult. We notice
that this Fuchs equation is an interesting alternative method to
prove the hyperbolicity of the limit cycle in system
(\ref{cicle4ch2}). This kind of equation may exist for all
algebraic limit cycle of a planar polynomial system and may let
distinguish its hyperbolic character.

\vspace{0.5cm}

\begin{center}
{\bf Acknowledgements}
\end{center}

The authors want to thank the fruitful conversations with
Professor L. Gavrilov from Univesit\'e de Toulouse which made
possible an easier proof for Theorem \ref{thhyp} as presented in
Section 3 of this work.


\begin{thebibliography}{99}

\bibitem{Andronov} {\sc A.A. Andronov, E.A. Leontovich, I.I.
Gordon and A.G. Ma\u{\i}er},{\it \  Qualitative theory of
second-order dynamic systems.} Translated from the Russian by D.
Louvish. Halsted Press (A division of John Wiley \& Sons), New
York-Toronto, Ont.; Israel Program for Scientific Translations,
Jerusalem-London, 1973.

\bibitem{Andronov2} {\sc A.A. Andronov, E.A. Leontovich, I.I.
Gordon and A.G. Ma\u{\i}er},{\it \  Theory of bifurcations of
dynamic systems on a plane.} Translated from the Russian. Halsted
Press (A division of John Wiley \& Sons), New York-Toronto, Ont.;
Israel Program for Scientific Translations, Jerusalem-London,
1973.

\bibitem{Arnold1} {\sc V.I. Arnol'd, S.M. Guse\u\i n-Zade and A.N. Varchenko},
{\it \  Singularities of differentiable maps. Vol. I. The
classification of critical points, caustics and wave fronts.}
Translated from the Russian by Ian Porteous and Mark Reynolds.
Monographs in Mathematics, {\bf 82}. Birkh\"{a}user Boston, Inc.,
Boston, MA, 1985.

\bibitem{Arnold2} {\sc V.I. Arnol'd, S.M. Guse\u\i n-Zade and A.N. Varchenko},
{\it \ Singularities of differentiable maps. Vol. II.
Monodromy and asymptotics of integrals.} Translated from the
Russian by Hugh Porteous. Translation revised by the authors and
James Montaldi. Monographs in Mathematics, {\bf 83}.
Birkh\"{a}user Boston, Inc., Boston, MA, 1988.

\bibitem{Chavarriga} {\sc J. Chavarriga},{\it \ A new example of quartic
algebraic limit cycle for a quadratic system}, Preprint,
Universitat de Lleida (1998).

\bibitem{ChGiLl} {\sc J. Chavarriga, H. Giacomini and J. Llibre}, {\it \ Uniqueness of algebraic limit cycles for quadratic
systems.}, Journal of Mathematical Analysis and Applications, {\bf
261}, (2001), 85--99.

\bibitem{chavallibre} {\sc J. Chavarriga and J. Llibre},{\it \ Invariant
algebraic curves and rational first integrals for planar
polynomial vector fields.} J. Differential Equations {\bf 169}
(2001), 1--16.

\bibitem{ChLlMo} {\sc J. Chavarriga, J. Llibre and J. Moulin-Ollagnier},{\it \ On a
result of Darboux.} LMS J. Comput. Math. {\bf 4} (2001), 197--210
(electronic).

\bibitem{ChLlSo} {\sc J. Chavarriga, J. Llibre, and J. Sorolla},
{\it \ Algebraic limit cycles of degree $4$ for quadratic
systems.} J. of Differential Equations {\bf 200} (2004), 206--244.


\bibitem{Chin} {\sc Yuan-sh\"un Ch'in}, {\it On algebraic limit cycles of degree $2$
of the differential equation $${dy\over dx}={\sum\sb{0\leqq
i+j\leqq 2}a\sb {ij}x\sp iy\sp j\over \sum\sb{0\leqq i+j\leqq
2}b\sb {ij}x\sp iy\sp j}.$$} Sci. Sinica {\bf 7} (1958), 934--945,
and Acta Math. Sinica {\bf 8} (1958), 23--35.

\bibitem{ChLlSw} {\sc C. Christopher, J. Llibre and G. \'{S}wirszcz},{\it
\ Invariant algebraic curves of large degree for quadratic
systems.} J. Math. Anal. Appl., in press.

\bibitem{Evdokimenko1}{\sc R.M. Evdokimenco}, {\it \ Construction
of algebraic paths and the qualitative investigation in the large
of the properties of integral curves of a system of differential
equations}, Differential Equations {\bf 6} (1970), 1349--1358.

\bibitem{Evdokimenko2}{\sc R.M. Evdokimenco}, {\it \ Behavior of
integral curves of a dynamic system}, Differential Equations {\bf
9} (1974), 1095--1103.

\bibitem{Evdokimenko3}{\sc R.M. Evdokimenco}, {\it \ Investigation in the
large of a dynamic system with a given integral curve},
Differential Equations {\bf 15} (1979), 215--221.

\bibitem{Filipstov}{\sc V.F. Filiptsov}, {\it \ Algebraic limit cycles},
Differential Equations {\bf 9} (1973), 983--988.

\bibitem{Perko} {\sc L. Perko},{\it \ Differential equations
and dynamical systems.} Third edition. Texts in Applied
Mathematics, {\bf 7}. Springer-Verlag, New York, 2001.

\bibitem{Yablonskii} {\sc A.I. Yablonskii},
{\it \ On limit cycles of certain differential equations},
Differential Equations {\bf 2} (1966), 164--168.

\bibitem{Ye} {\sc Ye Yan Qian, Cai Sui Lin, Chen Lan Sun,
Huang Ke Cheng, Luo Ding Jun, Ma Zhi En, Wang Er Nian, Wang Ming
Shu, Yang Xin An},{\it \ Theory of limit cycles.} Translated from
the Chinese by Chi Y. Lo. Second edition. Translations of
Mathematical Monographs, {\bf 66}. American Mathematical Society,
Providence, RI, 1986.


\end{thebibliography}
\end{document}